\documentclass[10pt]{article}

\usepackage[T1]{fontenc}
\usepackage[utf8]{inputenc}
\usepackage[a4paper,margin=1in]{geometry}
\usepackage[tbtags]{amsmath}
\usepackage{amssymb,amsthm,bm}
\usepackage{xcolor}
\usepackage{microtype}
\usepackage[colorlinks=true,linkcolor=blue,citecolor=blue,urlcolor=blue]{hyperref}

\allowdisplaybreaks[4]
\numberwithin{equation}{section}
\newcommand{\COM}[1]{}

\theoremstyle{plain}
\newtheorem{theo}{Theorem}[section]
\newtheorem{sat}[theo]{Proposition}
\newtheorem{lem}[theo]{Lemma}
\newtheorem{korr}[theo]{Corollary}
\theoremstyle{definition}
\newtheorem{de}[theo]{Definition}
\newtheorem{exxa}[theo]{Example}
\newtheorem{example}[theo]{Example}
\theoremstyle{remark}
\newtheorem{remark}[theo]{Remark}
\newtheorem{remarks}[theo]{Remarks}
 
\newcommand{\nelem}[1]{Lemma~\ref{#1}}

\newcommand{\netheo}[1]{Theorem~\ref{#1}}
\newcommand{\nekorr}[1]{Corollary~\ref{#1}}

\newcommand{\prooftheo}[1]{\noindent\textsc{Proof of Theorem}~\ref{#1}. }

\newcommand{\prooflem}[1]{\noindent\textsc{Proof of Lemma}~\ref{#1}. }
\newcommand{\proofkorr}[1]{\noindent\textsc{Proof of Corollary}~\ref{#1}. }

\newcommand{\kb}[1]{\bm{#1}}
\newcommand{\vk}[1]{\kb{#1}}
\newcommand{\kal}[1]{\mathcal{#1}}

\newcommand{\ve}{\varepsilon}
\newcommand{\fracl}[2]{\bigl(#1/#2\bigr)}

\newcommand{\abs}[1]{\lvert #1 \rvert}

\newcommand{\norm}[1]{\lVert #1 \rVert}
\newcommand{\normp}[1]{\lVert #1 \rVert}
\newcommand{\normS}[1]{\lVert #1 \rVert}

\newcommand{\E}[1]{\mathbf{E}\{#1\}}

\newcommand{\pk}[1]{\mathbf{P}\{#1\}}

\newcommand{\R}{\mathbb{R}}

\newcommand{\inr}{\in \R}

\newcommand{\ldot}{,\ldots,}
\newcommand{\limit}[1]{\lim_{#1\to\infty}}

\newcommand{\equaldis}{\stackrel{d}{=}}

\DeclareMathOperator*{\argmin}{argmin}

\newcommand{\BQN}{\begin{eqnarray}}
\newcommand{\EQN}{\end{eqnarray}}
\newcommand{\BQNY}{\begin{eqnarray*}}
\newcommand{\EQNY}{\end{eqnarray*}}
\newcommand{\BS}{\begin{sat}}
\newcommand{\ES}{\end{sat}}
\newcommand{\BT}{\begin{theo}}
\newcommand{\ET}{\end{theo}}
\newcommand{\BK}{\begin{korr}}
\newcommand{\EK}{\end{korr}}
\newcommand{\BD}{\begin{de}}
\newcommand{\ED}{\end{de}}
\newcommand{\BIT}{\begin{itemize}}
\newcommand{\EIT}{\end{itemize}}
\newcommand{\BDI}{\begin{description}}
\newcommand{\EDI}{\end{description}}
\newcommand{\BRM}{\begin{remarks}}
\newcommand{\ERM}{\end{remarks}}
\newcommand{\BTH}{\begin{theo}}
\newcommand{\ETH}{\end{theo}}
\newcommand{\BPR}{\begin{sat}}
\newcommand{\EPR}{\end{sat}}
\newcommand{\BEX}{\begin{exxa}}
\newcommand{\EEX}{\end{exxa}}
\newcommand{\BC}{\begin{cases}}
\newcommand{\EC}{\end{cases}}
\newcommand{\BL}{\begin{lem}}
\newcommand{\EL}{\end{lem}}

\newcommand{\QED}{\hfill $\Box$}

\newcommand{\IF}{\infty}

\newcommand{\SI}{\Sigma}
\newcommand{\SIJJ}{\Sigma_{JJ}}
\newcommand{\SIII}{\Sigma_{II}}

\newcommand{\SIJI}{\Sigma_{JI}}
\newcommand{\SIM}{\SI^{-1}}
\newcommand{\SIJJM}{\SIJJ^{-1}}
\newcommand{\SIIIM}{\Sigma_{II}^{-1}}

\newcommand{\njk}{\{1\ldot k\}}

\newcommand{\x}{\vk{x}}
\newcommand{\y}{\vk{y}}
\newcommand{\X}{\vk{X}}
\newcommand{\Y}{\vk{Y}}

\newcommand{\U}{\vk{\kal{U}}}

\newcommand{\IB}{I}
\newcommand{\JB}{J}
\def\1d{\{1\ldot d\}}

\newcommand{\qp}[2]{\kal{P}(#1,#2)}

\newcommand{\nouI}{\normS{\thrI}}

\newcommand{\thr}{\vk{u}_n}

\newcommand{\limuJ}{\vk{q}_J}
\newcommand{\limu}{\vk{q}}
\newcommand{\thry}{\thr}

\newcommand{\ntny}{\delta_n}
\newcommand{\wtny}{\zeta_n}
\newcommand{\wny}{\vk{v}_n}

\newcommand{\yJ}{\vk{y}_J}
\newcommand{\normyJ}{\yJ^\top(\SIM)_{JJ}\yJ}

\newcommand{\TSY}{\vk{u}^\top\SIM\y}

\newcommand{\vkal}{\vk{\alpha}}

\newcommand{\GSDAKF}{\kal{GSD}(k,\vkal,F)}

\newcommand{\LPSDA}{\kal{SD}(k,\vkal)}

\newcommand{\barAL}{\overline{\alpha}}
\newcommand{\barALI}{\overline{\alpha}_I}

\newcommand{\barALL}{\overline{\alpha}_L}

\newcommand{\bs}{\overline{\vk{a}}}
\newcommand{\bb}{\vk{a}}
\renewcommand{\b}{\bb}

\newcommand{\vkTI}{\vk{u}_I}

\title{Exact Tail Asymptotics of Dirichlet Distributions}
\author{Enkelejd Hashorva }

\begin{document}
\maketitle

\begin{abstract}
Let $\X=A^\top R\U$ be a linearly transformed generalised symmetrised Dirichlet scale mixture in $\R^k$, $k\ge2$.  For a fixed direction $\b\in(0,\infty)^k$, we derive an exact asymptotic expansion of $\pk{\X>\vk t_n}$ for eventually positive threshold vectors $\vk t_n$ described relative to the quadratic-programming minimiser on the natural active and residual Gumbel scales; residual limits equal to $-\infty$ are allowed.  The radial distribution is assumed to belong to the Gumbel max-domain of attraction.  The local power and constant are determined by the local product-power behaviour of the angular density near the minimising direction.  The result includes the ray $\vk t_n=u_n\b$ and yields an explicit comparison with the associated elliptical model, a conditional weak limit for the locally rescaled vector and the limiting location of the smallest component under a high common threshold.  The minimum overshoot is asymptotically exponential and independent of its location.  The finite-dimensional Gaussian minimum and location limits are recovered as a special case.
\end{abstract}


\section{Introduction}
Let $\X=(X_1\ldot X_k)^\top$ be a random vector in $\R^k$, $k\ge 2$, and let $u_n$, $n\ge 1$, be a positive sequence such that $u_n\to\infty$.   For any vector $\b=(a_1\ldot a_k)^\top\in(0,\infty)^k$, the events
\[
 \{\X>u_n\b\}:=\{X_1>u_n a_1\ldot X_k>u_n a_k\},\qquad n\ge 1,
\]
are tail events, that is
\[
 \limit{n}\pk{\X>u_n\b}=0.
\]
For such events it is natural to determine the exact rate at which $\pk{\X>u_n\b}$ converges to zero.

If $\X$ is a centred non-degenerate Gaussian random vector in $\R^k$ with covariance matrix $\SI$, then it admits the following stochastic representation; see, for example, \cite{cambanis}:
\begin{eqnarray}\label{eq:stochG}
 \X &\equaldis& R A^\top\U,
\end{eqnarray}
where $R>0$ is such that $R^2$ has a chi-square distribution with $k$ degrees of freedom, $A$ is a square matrix satisfying $A^\top A=\SI$, and $\U=(\kal{U}_1\ldot\kal{U}_k)^\top$ is uniformly distributed on the unit sphere of $\R^k$ and is independent of $R$. Tail asymptotics for Gaussian random vectors are well understood; see, e.g., \cite{BERMAN62,hh2003,EH2005a,EH19} and the references therein.

Hereafter we shall assume that the real matrix $A$ is non-singular. 

The radial decomposition in \eqref{eq:stochG} is also useful beyond the Gaussian case. If the distribution function $F$ of $R$ is left unspecified, then $\X$ is an elliptical random vector. If $F$ belongs to the Gumbel max-domain of attraction, see \eqref{eq:rdfd} below, then Theorem 3.1 of \cite{EH2007} implies
\begin{eqnarray}\label{eq:tailF}
 \pk{\X>u_n\b} &=& (1+o(1))\Psi(u_n)\pk{R>\mu u_n},\qquad n\to\infty,
\end{eqnarray}
for any $\b\in\R^k\setminus(-\infty,0]^k$, where $\Psi$ is an explicit function and
\[
 \mu=(\bs^\top\SIM\bs)^{1/2}.
\]
Here $\bs$ is the unique solution of the quadratic programming problem
\begin{eqnarray}\label{sol:x}
 \qp{\SIM}{\b}:\quad \text{minimise } \x^\top\SIM\x \text{ subject to } \x\in[a_1,\infty)\times\cdots\times[a_k,\infty).
\end{eqnarray}

The aim of this paper is to extend \eqref{eq:tailF} to random vectors with the same radial form as in \eqref{eq:stochG}, but with $\U$ following a symmetrised Dirichlet distribution with parameter $\vkal\in(0,\infty)^k$. We call $R\U$ a generalised symmetrised Dirichlet scale mixture and study its linear transform $A^\top R\U$; the underlying class was introduced by Fang and Fang \cite{FangFang1990}.

The main contribution is an exact asymptotic expansion of $\pk{\X>u_n\b}$ as $n\to\infty$ in this Dirichlet setting. The expansion has the same general structure as \eqref{eq:tailF}, but its local factor depends explicitly on $\vkal$.  The leading radial level $\mu u_n$ remains determined by the quadratic programming problem and does not depend on the vector $\vkal$.

The mechanism can already be read from the angular law.  With respect to Euclidean surface measure on the unit sphere, the symmetrised Dirichlet direction has density
\[
 h_{\vkal}(\vk u)=\frac{\Gamma(\barAL)}{2\prod_{i=1}^k\Gamma(\alpha_i)}
 \prod_{i=1}^k|u_i|^{2\alpha_i-1}.
\]
The Gumbel radial tail localises the event near the quadratic-programming direction, while this product-power density determines the angular mass on the active and residual local scales.  In particular, a zero coordinate at the minimising direction may change the power of the radial normalisation, and dependent singular factors may create logarithmic corrections.

An additional consequence is that a transformed Dirichlet scale mixture $\X$ and its associated elliptical random vector
\[
 \X^*\equaldis A^\top R\vk{V},
\]
where $\vk{V}$ is uniformly distributed on the unit sphere of $\R^k$ and is independent of $R$, may have the same tail asymptotic behaviour up to a multiplicative constant. This happens when the index set
\[
 \{i:(C\bs)_i=0,\ \alpha_i\ne 1/2,\ 1\le i\le k\},
\]
is empty, where $C=(A^\top)^{-1}$.  In that case the ratio of the two tail probabilities converges to the explicit constant in Corollary~\ref{korr:ell-comparison}.

The present paper should be read as a Dirichlet-scale-mixture analogue of the type I elliptical theory developed in \cite{EH2007}.  The Kotz Type III elliptical case, in which the radial tail has the form $p u^N\exp(-q u^\delta)$, was treated in \cite{EH2009Kotz}.  Theorem~\ref{theo:main1} below recovers the positive-direction case of that result when $\alpha_i=1/2$ for all $i$, and extends it to Dirichlet angular distributions.  The papers \cite{EH2007,EH2009Kotz} also used exact orthant tails for deriving approximations of conditional excesses.  Here we keep only the tail-ratio consequence that is genuinely affected by the Dirichlet angular law.  We do not pursue density approximation in this paper.

For the common threshold $\X>u_n\vk{1}$, the conditional rescaled-vector limit also identifies the smallest component and its location.  Under the hypotheses of the main theorem, the limiting minimum overshoot is exponential and independent of its location.  The limiting location is supported by the active set of the quadratic programme and is generally tilted by the product-power angular density.  When no active-scale angular zero is present, the location weights reduce to the usual quadratic-programming weights.  In particular, taking $\alpha_i=1/2$ and a $\chi_k$-distributed radius gives exactly the finite-dimensional Gaussian conditional-excess, minimum-overshoot and minimiser-location limits.

The paper is organised as follows. Section~2 contains notation and preliminary facts. The main result and three diagnostic examples are given in Section~3. Section~4 derives the conditional excess, rescaled-vector and minimum-location limits. The proofs of the main results are collected in Section~5, followed by an appendix.

\section{Preliminaries}
We first introduce some notation. Let $I$ be a non-empty subset of $\njk$ and put $J:=\njk\setminus I$, which may be empty.  For $\x=(x_1 \ldot x_k)^\top \inr^k$, define the subvector of $\x$ indexed by $I$ by
$\x_I:=(x_i, i \in I)^\top \inr^{\abs{I}}$. If $\SI\in \R^{k\times k}$ is
a square matrix, then the matrix  $\SI_{IJ}$ is obtained by retaining
both the rows and the columns of $\SI$ with indices in $I$ and in
$J$, respectively. Similarly we define $\SI_{JI}, \SI_{JJ},
\SI_{II}$. For notational simplicity we write $\x_I^\top, \SIJJM$
instead of $(\x_I)^\top, (\SI_{JJ})^{-1}$, respectively.
Given $\x,\y\inr^k$ we define
\BQNY
 \x&>&\y, \text{ if } x_i>y_i,\quad  \forall\, i=1 \ldot k,\\
\x&\ge& \y, \text{ if } x_i\ge y_i, \quad \forall\,  i=1 \ldot k,\\
\x+ \y &:=& (x_1+ y_1 \ldot x_k+ y_k)^\top,\\
c \x&:=&(cx_1\ldot cx_k)^\top,\quad c\inr,\\
\x\y&:=& ( x_1 y_1 \ldot x_k y_k)^\top, \quad
\x/\vk{y}:= (x_1/y_1 \ldot x_k/y_k)^\top.\\
\normS{\x_I}^2&:=& \x_I^\top \SIIIM\x_I.
  \EQNY
The coordinates of vectors in $\R^I$ are indexed by the elements of $I$; thus, for example, $(\SIIIM\x_I)_i$ is unambiguous for $i\in I$.

We write $\kal{B}_{a,b}$ for a Beta random variable with positive
parameters $a$ and $b$ and density
\[
\frac{\Gamma(a+b)}{\Gamma(a)\Gamma(b)}x^{a-1}(1-x)^{b-1},
\qquad x\in(0,1),
\]
where $\Gamma(\cdot)$ denotes the Gamma function.  Further, we write
$\vk{Y}\sim H$ if the random vector $\vk{Y}\inr^k$, $k\ge 1$, has
distribution function $H$.  For a univariate distribution function $H$,
set $\overline{H}:=1-H$.

Throughout this paper $\vkal:=(\alpha_1 \ldot \alpha_k)^\top $ stands for a vector in $\R^k$
with positive components, and
\[
\barAL:=\sum_{i=1}^k \alpha_i,\qquad
\overline{\alpha}_K:=\sum_{i\in K}\alpha_i,\qquad K\subset\njk.
\]
When $K$ is empty, $\overline{\alpha}_K$ equals $0$; throughout the paper, a product over an empty index set equals $1$.

\BD
Let $\vk D=(D_1\ldot D_k)^\top$ have the Dirichlet distribution with parameter $\vkal$, and let $\varepsilon_1\ldot\varepsilon_k$ be independent Rademacher random variables, independent of $\vk D$.  The random vector
\[
 \U=(\varepsilon_1\sqrt{D_1}\ldot\varepsilon_k\sqrt{D_k})^\top
\]
is said to have the symmetrised Dirichlet distribution with parameter $\vkal$; we write $\U\sim\LPSDA$.  In particular, $\U^\top\U=1$ almost surely.  This construction implies that $(\kal{U}_1 \ldot \kal{U}_{k-1})^\top$ has density
\BQN\label{eq:fang:1}
h(u_1 \ldot u_{k-1})&:=& \frac{\Gamma(\barAL)}{ \prod_{i=1}^k
\Gamma(\alpha_i)} \Bigl(1- \sum_{i=1}^{k-1} u_i^2 \Bigr)^{\alpha_k
-1}\prod_{i=1}^{k-1} \abs{u_i}^{ 2\alpha_i -1},   \quad
\sum_{i=1}^{k-1} u_i^2\le 1.
\EQN
\ED
For $k=1$, we use the convention that $\kal{SD}(1,\alpha_1)$ is the Rademacher distribution on $\{-1,1\}$.
Note that if $\vkal= \vk{1}/2$ with $ \vk{1}:=(1\ldot 1)^\top \inr^k$, then $\U \sim \LPSDA$ is uniformly distributed on the unit sphere of $\R^k$.

\BD  \label{def:LPSD} A random vector $\X_0$ in $\R^k$, $k\ge 1$, is called a generalised symmetrised Dirichlet scale mixture if
$\X_0\equaldis R\U$, where $R>0$ almost surely is independent of $\U$, $R\sim F$, and $\U\sim\LPSDA$. We write $\X_0\sim\GSDAKF$.  For a non-singular matrix $A$, the vector $A^\top\X_0$ is called a linearly transformed generalised symmetrised Dirichlet scale mixture.
 \ED
We henceforth consider distribution functions $F$ with an infinite upper endpoint.
Any subvector of a $\kal{GSD}$ scale mixture is again a $\kal{GSD}$ scale mixture; see \nelem{lem:dist} in Appendix.  A useful benchmark is obtained by taking
$R^2\sim\operatorname{Gamma}(\barAL,1/2)$.  Then $\X_0=R\U$ has independent components with
\BQN\label{eq:k:G}
 \abs{(X_0)_i}^2 \sim \operatorname{Gamma}(\alpha_i,1/2), \quad \forall i=1 \ldot k,
 \EQN
where $\operatorname{Gamma}(\alpha_i,1/2)$ denotes the Gamma distribution
with shape $\alpha_i$ and rate $1/2$.  Equivalently, $\X_0$ has density
\[
 h(\x)=\frac{2^{-\barAL}}{\prod_{i=1}^k\Gamma(\alpha_i)}
 \exp(-\x^\top\x/2)\prod_{i=1}^k|x_i|^{2\alpha_i-1},
 \qquad \x\in\R^k.
\]
For $\alpha_i=1/2$, $1\le i\le k$, this is the standard Gaussian density.

\section{Main Result}
Consider a linearly transformed $\kal{GSD}$ scale mixture $\X$ in $\R^k$, $k\ge 2$, with stochastic representation
\BQN \label{eq:def:gsd}
\X\equaldis  A^\top R  \U,
\EQN
where $R\sim F$ is independent of $\U \sim \kal{SD}(k, \vkal)$, $\vkal\in (0,\IF)^k$, and $A$
is a non-singular $k$-dimensional square matrix. Put $C:=(A^\top)^{-1}$. The positive definite matrix
$\Sigma:=A^\top A$ is the shape matrix of the model; except in the elliptical case, it need not be proportional to the covariance matrix of $\X$, when that covariance exists.\\
We begin with the threshold sequence $u_n\b$, $n\ge1$, and the joint
survivor probability $\pk{\X>u_n\b}$.  As in the elliptical setup
\cite{EH2007}, its tail asymptotics are closely related to the solution of
the quadratic programming problem \eqref{sol:x}.
If
 \BQN
 \label{eq:savage}
\SIM \b &> &\vk{0}, \quad  \vk{0}:=(0\ldot 0)^\top \inr^k
\EQN
is satisfied, then the minimum of the quadratic programming problem \eqref{sol:x} is attained at $\b$,
otherwise there exists a unique non-empty index set $I \subset \{1 \ldot k \}$
which defines the unique solution $\bs$ of $\qp{\SIM}{\b}$, see \nelem{prop:pre:gaus} in Appendix.
In the following we refer to  the index set $I$ as the minimal index set.\\
The only tail assumption imposed below on $F$ is that it belongs to the Gumbel max-domain of attraction with
some scaling function $w$ (for short $F \in GMDA(w)$), i.e.,
\BQN
\label{eq:rdfd}
\limit{u} \frac{\overline{F}(u+x/w(u))} {\overline{F}(u)} &=& \exp(-x),\quad \forall x\inr,
\EQN
with $\overline{F}:=1-F$. See \cite{falk,dehaan,Reiss1989,resnick2}
for standard accounts of max-domains of attraction.  For a recent
treatment of hidden regular variation and multivariate tail risk, see
\cite{resnick3}.

\noindent\textbf{Assumption A1.}
$F$ is a univariate distribution function with infinite upper endpoint
such that $F(0)=0$ and \eqref{eq:rdfd} holds.

We use without further mention that the convergence in \eqref{eq:rdfd}
is locally uniform in $x$ and that the reciprocal auxiliary function
$1/w$ is self-neglecting; equivalently,
\[
 \frac{w(u+x/w(u))}{w(u)}\longrightarrow1, \quad u \to \IF  
\]
locally uniformly for $x\in\R$.  Both facts are standard consequences
of \eqref{eq:rdfd}; see, e.g., \cite{dehaan,resnick2}.  We shall also
use two further standard consequences of Assumption A1.  First, we have that 
\BQN \label{eq:uv}
\limit{u} u w(u)&=& \IF.
\EQN
Second, for any $r>1$ and $\eta\in\R$
\BQN\label{eq:resn}
\limit{x} \frac{(x w(x))^\eta \overline{F}(rx)}{ \overline{F}( x)}&=&0.
\EQN
A proof of the latter rapid-variation bound is recalled in Appendix~\ref{sec:auxiliary}.

In the sequel $u_n$, $n\ge 1$, is a sequence of constants converging to infinity, and $\b\in(0,\infty)^k$ is a given vector.
If \eqref{eq:rdfd} holds, then we set
\BQN \label{main:def}
\ntny:= u_n \normS{\b_I},   \quad \wtny:= w(\ntny), \quad \lambda_n:= \ntny\wtny, \quad n\ge 1,
\EQN
with $I$ the minimal index set of $\qp{\SIM}{\bb}$.  By
\nelem{prop:pre:gaus}, $\normS{\b_I}\in(0,\IF)$.
The parameter $\vkal$ also plays a crucial role in these tail asymptotics,
through the following two index sets:
\BQN \label{LM}
L:=\{1 \le  i\le k: \alpha_i\not=1/2, \quad (C\bs)_i=0\}, \quad M:=\njk \setminus L
\EQN
appear explicitly in our asymptotic expansion. When $L$ is empty, $\X$ and the associated elliptical random vector $\X^*\equaldis A^\top R \vk{V}$, with $\vk{V}\sim \kal{SD}(k, \vk{1}/2)$ independent of $R$, have the same tail asymptotics up to the explicit constant in Corollary~\ref{korr:ell-comparison}. We state next our main result.

\def\nouI{\ntny}
\def\thrz{\thry}

\BT\label{theo:main1}
Let $\X$ be the linearly transformed $\kal{GSD}$ scale mixture in $\R^k$, $k\ge 2$, defined in \eqref{eq:def:gsd}, and set
$C:=(A^\top)^{-1}, \SI:= A^\top A$. For a
given vector $\bb\in(0,\infty)^k$, let $I$ with $m$ elements be the minimal index set corresponding to $\qp{\SIM}{\bb}$ with the unique solution $\bs\ge\bb$.  Put $J:=\njk\setminus I$. Suppose that Assumption A1 holds, and define $\ntny, \wtny, \lambda_n$ as in \eqref{main:def}.
Let $\thry$, $n\ge 1$, be a sequence of thresholds in $\R^k$ that is componentwise positive for all sufficiently large $n$ and satisfies
\BQN \label{eq:M}
\limit{n}\wtny \bigl( \thry-  u_n \bs\bigr)_I &=& \vk{q}_I \inr^{m},
\EQN
and if $m< k$
\BQN\label{eq:wtn:thr}
\limit{n}\fracl{\wtny}{\ntny }^{1/2} \bigl(\thry  - u_n\bs \bigr)_J
&=& \limuJ\in [-\infty,\IF)^{k-m}.
\EQN
If $J\ne\emptyset$, let $\vk q$ be the unique vector with subvectors $\vk q_I$ and $\vk q_J$, with the evident interpretation of components equal to $-\infty$; if $J=\emptyset$, put $\vk q=\vk q_I$. All terms below involving an empty $J$ are omitted. For $K\subset\njk$ write $C_{iK}:=(C_{ij},j\in K)$ for the row vector obtained from the $i$th row of $C$, and put
\BQN\label{eq:Lsplit}
L_J^0:=\{i\in L: J\ne\emptyset,\ C_{iJ}\ne \vk{0}_J\},\quad
L_I^0:=L\setminus L_J^0 .
\EQN
Put
\[
 L_J^-:=\{i\in L_J^0:\alpha_i<1/2\},\qquad
 L_I^-:=\{i\in L_I^0:\alpha_i<1/2\}.
\]
Assume that the matrices $C_{L_J^-,J}$ and $C_{L_I^-,I}$ have full row rank; the corresponding condition is void when the row set is empty.  Define
\BQN\label{eq:kappa}
\kappa:=1-\abs{I}-\frac{\abs{J}}{2}+\frac{\abs{L_J^0}}{2}+\abs{L_I^0}
       -\overline{\alpha}_{L_J^0}-2\overline{\alpha}_{L_I^0},
\EQN
and set, with $\vk{u}:=\b/\normS{\b_I}$
\BQNY
\tau_{L_I^0,L_J^0}&:=&\int_{\y_I>\vk{q}_I,\,\y_J>\limuJ}
\prod_{i\in L_I^0}\abs{C_{iI}\y_I}^{2\alpha_i-1}
\prod_{i\in L_J^0}\abs{C_{iJ}\y_J}^{2\alpha_i-1} \\
&&\hspace*{2.0cm}\times
\exp\left(-\vk{u}_I^\top\SIIIM\y_I-\frac12\y_J^\top(\SIM)_{JJ}\y_J\right)\,d\y_I\,d\y_J,\\
\tau_M^*&:=&\prod_{\substack{i\in M\\ \alpha_i\ne 1/2}}
\left(\normS{\b_I}^{1-2\alpha_i}\abs{(C\bs)_i}^{2\alpha_i-1}\right).
\EQNY
Then $\tau_{L_I^0,L_J^0}\in(0,\IF)$, and as $n\to\IF$
\BQN\label{eq:theo1:gen}
\pk{\X > \thry}
&\sim & \tau_{L_I^0,L_J^0}\tau_M^*
\frac{\Gamma(\barAL)}{2\prod_{i=1}^k\Gamma(\alpha_i)\abs{\SI}^{1/2}}
\lambda_n^{\kappa}\overline F(\ntny).
\EQN
In particular, if $J\ne\emptyset$ and $L_I^0=\emptyset$, then \eqref{eq:theo1:gen} becomes
\BQN
\label{eq:theo1:1}
\pk{\X > \thry}
&\sim &\tau_{J,L}\tau_M^* \frac{\Gamma(\barAL) \exp(- \vk{q}_I^\top \SIIIM \vkTI)}{2\prod_{i=1}^k\Gamma(\alpha_i)  \abs{\SI}^{1/2} \prod_{i\in I} (\SIIIM \vkTI)_i}
\notag\\
&& \times
\lambda_n^{1-\abs{I}-\abs{J}/2+\abs{L}/2-\barALL}  \overline{F}(\ntny), \quad n\to \IF,
\EQN
where
\BQNY
 \tau_{J,L}&:=& \int_{\y_J> \limuJ} \prod_{i\in L} \abs{C_{iJ}\y_J}^{2 \alpha_i-1}\exp(-\normyJ/2) \,  d\y_J .
\EQNY
If $J$ is empty, then $L_I^0=L$, and \eqref{eq:theo1:gen} reads
\BQN
\label{eq:theo1:2}
\pk{\X > \thry}
&\sim &\tau_L \tau_M^*\frac{\Gamma(\barAL)}{2\prod_{i=1}^k\Gamma(\alpha_i)  \abs{\SI}^{1/2} }
 \lambda_n^{1-k+ \abs{L}-2\barALL} \overline{F}(\ntny),
\EQN
where $\tau_L:=\int_{\y>\limu } \prod_{i\in L} \abs{(C\y)_i}^{2 \alpha_i-1}\exp(-\TSY) \,  d\y \in (0,\IF)$.
\ET
\medskip
For the ray $\thry=u_n\b$, condition \eqref{eq:M} holds with $\vk q_I=\vk0_I$.  For $j\in J$, the corresponding residual limit is $q_j=0$ when $\bs_j=b_j$ and $q_j=-\infty$ when $\bs_j>b_j$.  Thus \netheo{theo:main1} includes the fixed-direction orthant tail announced in the introduction.
\par\medskip
\noindent What matters for a zero component $(C\bs)_i$ is whether the row $C_{iJ}$ is non-zero; if it is, the dominant perturbation is of $J$-scale, otherwise it is of $I$-scale.  The rank condition is imposed only on singular factors, namely those with $\alpha_i<1/2$; dependent zero factors with $\alpha_i>1/2$ cause no integrability problem.  If, after ordering the active coordinates first, $A$ is block upper triangular (in particular, if $A$ is upper triangular in that ordering), equivalently $C_{IJ}=0$, then $L_I^0=L\cap I$ and $L_J^0=L\cap J$, and the rank condition is automatic because $C_{II}$ and $C_{JJ}$ are nonsingular.  Without the negative-power rank condition, several singular hyperplanes may coalesce; then the integral above can be infinite and the lower-order local scale can change the power of $\lambda_n$ or introduce logarithmic factors.
\par\medskip

\BK\label{korr:ell-comparison}
Under the assumptions of \netheo{theo:main1}, let $\X^*:=A^\top R\vk V$, where $\vk V$ is uniform on the unit sphere of $\R^k$ and independent of $R$.  If $L=\emptyset$, then
\[\limit{n}
\frac{\pk{\X>\thry}}{\pk{\X^*>\thry}}
 =
 c_{\vkal,A,\b}:=
 \tau_M^*\frac{\Gamma(\barAL)\Gamma(1/2)^k}
 {\Gamma(k/2)\prod_{i=1}^k\Gamma(\alpha_i)}
 \in(0,\infty).
\]
\EK

\begin{remark}[Elliptical and Kotz Type III special cases]\label{rem:ell-kotz}
If $\alpha_i=1/2$ for all $i$, then $L=\emptyset$, $\barAL=k/2$ and $\tau_M^*=1$.  The local integral in \eqref{eq:theo1:gen} is then Gaussian, and \netheo{theo:main1} reduces to the positive-direction, eventually positive-threshold case of the unbounded-endpoint type I elliptical expansion in \cite[Theorem~3.1]{EH2007}.  If, in addition,
\[
 \overline F(u)=(1+o(1))p u^N\exp(-q u^\delta),\qquad
 p,q,\delta>0,\quad N\in\R,
\]
then $F\in GMDA(w)$ with $w(u)=(1+o(1))q\delta u^{\delta-1}$.  Substitution of this radial tail into the preceding elliptical special case gives the corresponding positive-direction case of the Kotz Type III elliptical expansion in \cite[Theorem~3.1]{EH2009Kotz}, after the notational change from the present moving threshold $\thry$ to $t_n\vk a+\x/\vk v_n$.  Thus the genuinely new feature of the present theorem is the Dirichlet angular part, encoded by $L_I^0,L_J^0$ and the local integral \eqref{eq:theo1:gen}, rather than a new radial tail class.
\end{remark}

\begin{remark}[Gaussian special case]\label{rem:gaussian}
Let $\X\sim N_k(\vk 0,\Sigma)$ with
$\Sigma=A^\top A$, where $A$ is a non-singular $k\times k$ real matrix.  In the representation \eqref{eq:def:gsd} this corresponds to
$\alpha_i=1/2$, $1\le i\le k$, and $R^2\sim\chi_k^2$.  Hence as $u\to\infty$ with $w(u)=u$ we have 
\[
 \overline F(u)=(1+o(1))\frac{2^{1-k/2}}{\Gamma(k/2)}
 u^{k-2}\exp(-u^2/2)  .
\]
Let $I,J,\bs$ be as in \netheo{theo:main1}, put $m:=|I|$
\[
 \mu:=\normS{\b_I}=(\b_I^\top\SIIIM\b_I)^{1/2},\qquad
 h_n:=u_n\mu,\qquad \vk u_I:=\b_I/\mu .
\]
Assume, as in \netheo{theo:main1}, that $\thry$ is eventually componentwise positive and that
\[
 h_n(\thry-u_n\bs)_I\to \vk q_I\in\R^m,
\]
and, if $J\ne\emptyset$, that
\[
 (\thry-u_n\bs)_J\to \vk q_J\in[-\infty,\infty)^{|J|}.
\] Then, as $n\to\infty$
\[
\pk{\X>\thry}
=(1+o(1))
\frac{
 \exp\{-h_n^2/2-\vk q_I^\top\SIIIM\vk u_I\}
 \pk{\X_J>\vk q_J\mid \X_I=\vk 0_I}
}
{(2\pi)^{m/2}\abs{\SIII}^{1/2}h_n^m
 \prod_{i\in I}(\SIIIM\vk u_I)_i},
\]
where the conditional probability is interpreted as one when $J=\emptyset$,
and lower limits equal to $-\infty$ are understood componentwise. 
\end{remark}
\par\medskip

\begin{example}[An independent-component benchmark]\label{ex:factorising}
Let $A=I_k$ and let $R^2\sim\operatorname{Gamma}(\barAL,1/2)$.  By
\eqref{eq:k:G}, the components of $\X=R\U$ are independent and
$|X_i|^2\sim\operatorname{Gamma}(\alpha_i,1/2)$.  Hence, for
$\bb=(a_1\ldot a_k)^\top\in(0,\infty)^k$,
\[
 \pk{X_i>u a_i}
 \sim \frac{u^{2\alpha_i-2}a_i^{2\alpha_i-2}}
 {2^{\alpha_i}\Gamma(\alpha_i)}
 \exp(-u^2a_i^2/2),
 \qquad u\to\infty,
\]
and independence implies as $n\to \IF$
\BQN\label{eq:factorising-direct}
 \pk{\X>u\bb}
 &\sim&
 \frac{u^{2\barAL-2k}}{2^{\barAL}\prod_{i=1}^k\Gamma(\alpha_i)}
 \left(\prod_{i=1}^k a_i^{2\alpha_i-2}\right)
 \exp\left(-\frac{u^2\norm{\bb}^2}{2}\right).
\EQN

This also checks every constant in \netheo{theo:main1}.  For this radial
law one may take $w(r)=r$.  Taking $u_n=u$, the quadratic programme has
$\bs=\bb$, $I=\njk$, $J=L=\emptyset$, and
\[
 \ntny=u\norm{\bb},\qquad \lambda_n=u^2\norm{\bb}^2.
\]
Moreover,
\[
 \tau_{\emptyset}
 =\int_{\y>\vk0}\exp\left(-\frac{\bb^\top\y}{\norm{\bb}}\right)d\y
 =\frac{\norm{\bb}^k}{\prod_{i=1}^k a_i},
 \qquad
 \tau_M^*=\norm{\bb}^{k-2\barAL}
 \prod_{i=1}^k a_i^{2\alpha_i-1},
\]
while
\[
 \overline F(r)\sim
 \frac{2^{1-\barAL}}{\Gamma(\barAL)}
 r^{2\barAL-2}e^{-r^2/2}.
\]
Substitution in \eqref{eq:theo1:2} yields exactly
\eqref{eq:factorising-direct}.  When all $\alpha_i=1/2$, this reduces to
the product of the usual Gaussian Mills ratios.
\end{example}
\par\medskip

\begin{example}[A two-dimensional transition]\label{ex:two-dimensional}
Let $k=2$ and
\[
 A=\begin{pmatrix}1&\rho\\0&d\end{pmatrix},\qquad
 d:=\sqrt{1-\rho^2},\qquad 0<\rho<1,
\]
so that
\[
 \Sigma=\begin{pmatrix}1&\rho\\ \rho&1\end{pmatrix},\qquad
 C=\begin{pmatrix}1&0\\-\rho/d&1/d\end{pmatrix},
\]
and
\[
 X_1=RU_1,\qquad X_2=R(\rho U_1+dU_2).
\]
Take $\bb=(1,c)^\top$ with $0<c<\rho$ and
$\vkal=(1/2,\gamma)^\top$, $\gamma>0$, $\gamma\ne1/2$.  The quadratic programme has
\[
 \bs=(1,\rho)^\top,\qquad I=\{1\},\qquad J=\{2\},
 \qquad C\bs=(1,0)^\top.
\]
Thus $L=L_J^0=\{2\}$ and $L_I^0=\emptyset$.  If $\gamma<1/2$, the
required rank condition holds since $C_{2J}=1/d\ne0$; if $\gamma>1/2$,
there is no negative-power rank condition.  Moreover, $\ntny=u_n$ and
$\lambda_n=u_n w(u_n)$.  For the fixed ray
$\thry=u_n\bb$ one has $q_I=0$ and $q_J=-\infty$.  Since
\[
 \int_0^\infty e^{-y_1}\,dy_1\,
 \times
 \int_{\R}\abs{y_2/d}^{2\gamma-1}
       \exp\left(-\frac{y_2^2}{2d^2}\right)dy_2
 =d\,2^\gamma\Gamma(\gamma)
\]
\netheo{theo:main1} implies as $n\to \IF$
\BQN\label{eq:two-dimensional-ray}
 \pk{\X>u_n\bb}
 &\sim&
 \frac{2^{\gamma-1}\Gamma(\gamma+1/2)}{\sqrt\pi}
 \lambda_n^{-\gamma}\overline F(u_n).
\EQN
This also follows directly: the second constraint is asymptotically slack and
$U_1^2\sim\kal{B}_{1/2,\gamma}$ while $\pk{U_1>0}=1/2$, so the right-hand side of
\eqref{eq:two-dimensional-ray} is the endpoint expansion of
$\pk{RU_1>u_n}$.

The moving thresholds in \netheo{theo:main1} describe the transition to the
boundary direction.  For fixed $q\in\R$, put
\[
 \vk t_n:=\left(u_n,\rho u_n+q\sqrt{\frac{u_n}{w(u_n)}}\right)^\top.
\]
Then the residual limit is $q_J=q$, and
\BQN\label{eq:two-dimensional-moving}
 \pk{\X>\vk t_n}
 &\sim&
 \frac{\Gamma(\gamma+1/2)}{2\sqrt\pi\,\Gamma(\gamma)}
 \left[\int_{q/d}^\infty
       \abs{z}^{2\gamma-1}e^{-z^2/2}\,dz\right]
 \lambda_n^{-\gamma}\overline F(u_n).
\EQN
For the fixed ray $u_n\bb$ with $c<\rho$, the residual limit is
$q_J=-\infty$, so its constant is obtained from
\eqref{eq:two-dimensional-moving} by replacing the lower limit by
$-\infty$.  If instead $\bb=(1,\rho)^\top$, its fixed ray has $q_J=0$ and
gives one half of the constant in \eqref{eq:two-dimensional-ray}.

The fixed-slope regime in \cite[Example~5]{EH2012Biv} is complementary.
For $\rho<a\le1$, put
\[
 \mu_a:=\frac{\sqrt{1-2a\rho+a^2}}{d},\qquad
 h_\gamma(s):=\frac{2(1-s^2)^{\gamma-1}}
                    {\mathrm B(1/2,\gamma)},\quad 0<s<1,
\]
where $h_\gamma$ is the density of $W:=|\kal U_1|$.  Specialising
\netheo{theo:main1} to $\bb=(1,a)^\top$ gives
\[
 \pk{X_1>u_n,\ X_2>a u_n}
 \sim
 \frac{d^2h_\gamma(1/\mu_a)}{4(1-a\rho)}
 \frac{\overline F(\mu_a u_n)}{u_nw(\mu_a u_n)}, 
 \qquad n\to\infty.
\]
This is exactly the fixed-$a>\rho$ expansion in
\cite[Example~5]{EH2012Biv}.  Thus that result covers the supercritical
fixed slopes, whereas \eqref{eq:two-dimensional-ray} treats the slack
regime $c<\rho$ and \eqref{eq:two-dimensional-moving} resolves the
boundary transition at $c=\rho$.

The conditional limit is equally explicit.  For $x_1,x_2\ge0$,
Corollary~\ref{korr:excess} gives
\begin{align}
 &\mathbf P\left\{
 w(u_n)(X_1-u_n)>x_1,\ 
 \sqrt{\frac{w(u_n)}{u_n}}\,(X_2-t_{n,2})>x_2
 \,\middle|\,\X>\vk t_n\right\}\notag \longrightarrow
 e^{-x_1}
 \frac{\displaystyle\int_{(q+x_2)/d}^{\infty}
 |z|^{2\gamma-1}e^{-z^2/2}\,dz}
 {\displaystyle\int_{q/d}^{\infty}
 |z|^{2\gamma-1}e^{-z^2/2}\,dz}
 \label{eq:two-dimensional-excess}
\end{align}
as $n\to \IF$. 
Thus the active excess is standard exponential and is independent of
the residual excess.  The latter is a polynomially tilted truncated
Gaussian law; in the associated elliptical case $\gamma=1/2$, the
polynomial tilt disappears.

This example also separates the roles of $\Sigma$ and $A$.  If $O$ is
orthogonal and $A$ is replaced by $OA$, one obtains a different
non-elliptical model with the same $\Sigma$ and the same quadratic programme,
but $C\bs$ is replaced by $OC\bs=O(1,0)^\top$.  For a generic $O$, the
second component of $O(1,0)^\top$ is non-zero, so $L=\emptyset$ and the
power is $\lambda_n^{-1/2}$ instead of $\lambda_n^{-\gamma}$.  There is no
contradiction: a symmetrised Dirichlet direction is not rotationally invariant
unless all parameters equal $1/2$.
\end{example}

\begin{example}[Coalescing angular hyperplanes]\label{ex:three-dimensional}
Failure of the negative-power rank condition in \netheo{theo:main1} can produce a genuine
lower-order effect.  Let $k=3$ and set
\[
 C=\begin{pmatrix}
 1&-2&1\\
 -2&1&1\\
 1&1&0
 \end{pmatrix},\qquad
 A^\top=C^{-1}=\begin{pmatrix}
 1/6&-1/6&1/2\\
 -1/6&1/6&1/2\\
 1/2&1/2&1/2
 \end{pmatrix}.
\]
Take
\[
 \bb=(1,1,1/2)^\top,\qquad
 \bs=(1,1,1)^\top,
 \qquad \vkal=(1/4,1/4,1/2)^\top.
\]
Since
\[
 C\bs=(0,0,2)^\top,\qquad
 \Sigma^{-1}\bs=C^\top C\bs=(2,2,0)^\top,
\]
the quadratic programme has minimal index set $I=\{1,2\}$ and $J=\{3\}$,
with $\ntny=2u_n$ and $\lambda_n=2u_n w(2u_n)$.  For the fixed ray,
$q_I=\vk0_I$ and $q_3=-\infty$, so both sides of the limiting hyperplane
$y_3=0$ contribute.  Moreover,
\[
 L=L_J^0=L_J^-=\{1,2\},\qquad
 C_{L_J^-,J}=\begin{pmatrix}1\\1\end{pmatrix},
\]
which has rank one.  The two limiting factors therefore coalesce into
\[
 \abs{y_3}^{-1/2}\abs{y_3}^{-1/2}=\abs{y_3}^{-1},
\]
and the putative local integral is logarithmically divergent.

The finite-$n$ factors show how the divergence is regularised.  With
$\varepsilon_n:=\lambda_n^{-1/2}$ they are
\[
 \abs{y_3+\varepsilon_n(y_1-2y_2)}^{-1/2}
 \abs{y_3+\varepsilon_n(-2y_1+y_2)}^{-1/2}.
\]
Their roots are separated by
$3\varepsilon_n\abs{y_1-y_2}$.  For fixed $y_1,y_2>0$, $y_1\ne y_2$,
put
\[
 s:=\frac{y_1+y_2}{2},\qquad r:=\frac{y_1-y_2}{2},\qquad
 z:=y_3-\varepsilon_n s.
\]
The angular product then becomes
$\abs{z+3\varepsilon_n r}^{-1/2}\abs{z-3\varepsilon_n r}^{-1/2}$, and
\[
 \int_{\R}
 \abs{y_3+\varepsilon_n(y_1-2y_2)}^{-1/2}
 \abs{y_3+\varepsilon_n(-2y_1+y_2)}^{-1/2}
 e^{-y_3^2}\,dy_3
 =\log\lambda_n+O\left(1+\abs{\log\abs{y_1-y_2}}\right).
\]
On every bounded box, the corresponding $z$-integral is bounded by
\[
 C\left(1+\log\frac1{\varepsilon_n}+\abs{\log\abs{r}}\right).
\]
Since $\abs{\log\abs{r}}$ is locally integrable, dominated convergence after division
by $2\log(1/\varepsilon_n)=\log\lambda_n$ gives the local logarithmic
coefficient.  The angular product has integral $O(\log\lambda_n)$ on every
unit cube, uniformly in its location.  Combining this estimate with the
arbitrary-order Potter bound \eqref{eq:radial-global-bound} makes the
complement of a growing ball $o(\log\lambda_n)$; the local convergence in
\eqref{eq:innerlimit} handles the ball itself.
The formal exponent given by \eqref{eq:kappa} is $\kappa=-1$,
$\abs{\Sigma}^{1/2}=1/6$, and
$\int_{(0,\infty)^2}e^{-y_1-y_2}\,dy_1dy_2=1$.  The preceding bounds justify
a direct modification of the localisation argument in the proof of
\netheo{theo:main1}, and yield
\BQN\label{eq:three-dimensional-log}
 \pk{\X>u_n\bb}
 &\sim&
 \frac{3}{\sqrt\pi\,\Gamma(1/4)^2}
 \lambda_n^{-1}\log\lambda_n\,\overline F(2u_n).
\EQN
Thus the negative-power rank hypothesis excludes a real logarithmic correction,
rather than merely simplifying the proof.  It is sufficient but not necessary: if
$\alpha_1+\alpha_2>1/2$, the corresponding coalesced power is locally
integrable.  Rows with $\alpha_i>1/2$ give zeros rather than singularities,
and no rank condition is imposed on them.  A sharper condition covering
integrable dependent singularities is possible, but would require a separate
uniform analysis of moving hyperplanes; we retain the transparent sufficient
condition in \netheo{theo:main1}.
\end{example}
\par\medskip

\section{Applications}

The applications below are obtained by taking ratios of the exact tail expansion in \netheo{theo:main1}.  This is the same mechanism used in the type I elliptical case in \cite{EH2007}; the difference here is that the Gaussian local limit is replaced by the Dirichlet-dependent local integral appearing in \eqref{eq:theo1:gen}.  Accordingly, all applications below follow from survivor-probability ratios.

For a threshold vector $\thry$ define the excess random vector
\[
 \X[\thry]\equaldis \X-\thry\mid \X>\thry .
\]
The multiplication by $\wny$ below is componentwise, with
\[
 (\wny)_I=\wtny\vk{1}_I,\qquad
 (\wny)_J=\left(\frac{\wtny}{\ntny}\right)^{1/2}\vk{1}_J,
\]
where the second relation is omitted if $J=\emptyset$.

\BK\label{korr:excess}
Assume the conditions of \netheo{theo:main1} and, when $J\ne\emptyset$, suppose that $\vk q_J\in\R^{|J|}$.  For a vector $\vk r$ with subvectors $\vk r_I$ and $\vk r_J$, define
\BQNY
\mathcal T(\vk r)&:=&\int_{\y_I>\vk r_I,\,\y_J>\vk r_J}
\prod_{i\in L_I^0}\abs{C_{iI}\y_I}^{2\alpha_i-1}
\prod_{i\in L_J^0}\abs{C_{iJ}\y_J}^{2\alpha_i-1}\\
&&\hspace*{1.7cm}\times
\exp\left(-\vk{u}_I^\top\SIIIM\y_I-\frac12\y_J^\top(\SIM)_{JJ}\y_J\right)\,d\y_I\,d\y_J,
\EQNY
with all $J$-terms omitted when $J=\emptyset$. Then, for all $\x\in[0,\infty)^k$ at continuity points of the limiting survivor function,
\BQN\label{eq:excesslimit}
 \pk{\wny\X[\thry]>\x}
 &\longrightarrow&
 \frac{\mathcal T(\vk q+\x)}{\mathcal T(\vk q)},\qquad n\to\infty,
\EQN

If $L_I^0=\emptyset$, then the limit factorises as
\BQN\label{eq:excessfactor}
\frac{\mathcal T(\vk q+\x)}{\mathcal T(\vk q)}
&=&\exp(-\x_I^\top\SIIIM\vk{u}_I)
\frac{\displaystyle\int_{\y_J>\vk q_J+\x_J}\prod_{i\in L}\abs{C_{iJ}\y_J}^{2\alpha_i-1}\exp(-\y_J^\top(\SIM)_{JJ}\y_J/2)\,d\y_J}
{\displaystyle\int_{\y_J>\vk q_J}\prod_{i\in L}\abs{C_{iJ}\y_J}^{2\alpha_i-1}\exp(-\y_J^\top(\SIM)_{JJ}\y_J/2)\,d\y_J}.
\EQN
Consequently, the active part has independent exponential components with rates $(\SIIIM\vk u_I)_i$, $i\in I$, and is independent of the residual $J$-part.
\EK

If $(q_J)_j=-\infty$ for some $j\in J$, extend $\mathcal T$ by interpreting a lower bound equal to $-\infty$ as no restriction in that coordinate.  The same ratio argument then gives \eqref{eq:excesslimit} for finite $\x\ge\vk0$.  Its limiting survivor does not depend on $x_j$.  Equivalently, that coordinate of $\wny\X[\thry]$ diverges conditionally to $+\infty$ in the extended half-line; the law of the remaining finite coordinates is obtained by marginalising over all such coordinates.

Formula \eqref{eq:excesslimit} is the form used for conditional excess approximation: the numerator and denominator are the same local tail integral with shifted lower limits.  In particular, the Dirichlet contribution is not an additional radial correction; it changes the limiting conditional law through the factors involving the zero components of $C\bs$.

\begin{korr}[Conditional rescaled-vector and minimum-location limits]
\label{korr:minimum-location}
Assume the conditions of \netheo{theo:main1}, take
$\b=\vk{1}$ and $\thry=u_n\vk{1}$, and let $\bs$ and $I$ be the
quadratic-programming solution and its minimal index set.  Put
\[
 J:=\njk\setminus I,\qquad
 P:=\{j\in J:\bs_j=1\},\qquad
 N:=\{j\in J:\bs_j>1\},
\]
and
\[
 \mu^2:=\vk{1}_I^\top\SIIIM\vk{1}_I,\qquad
 \vk{\theta}:=\SIIIM\vk{1}_I,\qquad
 \sigma_*^2:=\mu^{-2}.
\]
Thus $\vk{\theta}>\vk{0}_I$.  Moreover, the vector
\[
 \vk{p}_I:=\frac{\vk{\theta}}{\mu^2},\qquad
 \vk{p}_J:=\vk{0}_J
\]
is the unique minimiser of
\[
 \min_{\vk{p}\ge\vk{0},\,\vk{1}^\top\vk{p}=1}
       \vk{p}^\top\Sigma\vk{p},
\]
and the minimum equals $\sigma_*^2$. Set
\[
 r_n:=\mu u_n,\qquad
 c_n:=\frac{w(r_n)}{\mu},\qquad
 d_n:=\left(\frac{w(r_n)}{r_n}\right)^{1/2}.
\]
For $\vk z_I\in(0,\infty)^{|I|}$ define
\[
 \psi_I(\vk z_I):=
 \prod_{\ell\in L_I^0}
       \abs{C_{\ell I}\vk z_I}^{2\alpha_\ell-1},
\qquad
 c_I:=\int_{\vk z_I>\vk0_I}
       \psi_I(\vk z_I)e^{-\vk\theta^\top\vk z_I}\,d\vk z_I.
\]
Let $\vk Z_I$ have density
\begin{equation}
\label{eq:active-limit-density}
 f_I(\vk z_I)=c_I^{-1}\psi_I(\vk z_I)
 e^{-\vk\theta^\top\vk z_I}
 \boldsymbol 1_{\{\vk z_I>\vk0_I\}}.
\end{equation}
If $J\ne\emptyset$, let $\vk Z_J^{(P)}$ have density proportional to
\begin{equation}
\label{eq:residual-limit-density}
 \boldsymbol 1_{\{\vk z_P>\vk0_P\}}
 \prod_{\ell\in L_J^0}
       \abs{C_{\ell J}\vk z_J}^{2\alpha_\ell-1}
 \exp\left\{-\frac12\vk z_J^\top(\SIM)_{JJ}\vk z_J\right\},
\end{equation}
where $\vk z_N$ ranges over $\R^{|N|}$; the positivity restriction is
void when $P=\emptyset$.  The normalising constants in
\eqref{eq:active-limit-density}--\eqref{eq:residual-limit-density}
are finite and positive.  Take $\vk Z_I$ and $\vk Z_J^{(P)}$ independent.
Then as $n\to\infty$
\begin{equation}
\label{eq:conditional-local-min}
 \mathcal L\left(
 \left(
 c_n(\X_I-u_n\vk{1}_I),
 d_n(\X_J-u_n\bs_J)
 \right)
 \,\middle|\,\min_{1\le i\le k}X_i>u_n\right)
 \Longrightarrow
 \mathcal L(\vk Z_I,\vk Z_J^{(P)}),
\end{equation}
with the $J$-block omitted when $J=\emptyset$.  In particular, we have the conditional limit as $n\to\infty$
\begin{equation}
\label{eq:conditional-active-passive}
 \mathcal L\left(
 \left(
 c_n(\X_I-u_n\vk{1}_I),
 d_n(\X_P-u_n\vk{1}_P),
 \frac{\X_N}{u_n}
 \right)
 \,\middle|\,\min_{1\le i\le k}X_i>u_n\right)
 \Longrightarrow
 \mathcal L(\vk Z_I,(\vk Z_J^{(P)})_P,\bs_N),
\end{equation}
where empty blocks are omitted.

Let
\[
 \underline X:=\min_{1\le i\le k}X_i,\qquad
 T:=\min\argmin_{1\le i\le k}X_i,
\]
and put
\[
 S:=\min_{i\in I}(Z_I)_i,\qquad
 T_*:=\argmin_{i\in I}(Z_I)_i.
\]
The argmin sets defining $T$ and $T_*$ are almost surely singletons and we have as $n\to \IF$ 
\begin{equation}
\label{eq:minimum-location}
 \mathcal L\left(
 \left(c_n(\underline X-u_n),T\right)
 \,\middle|\,\underline X>u_n\right)
 \Longrightarrow
 \mathcal L(S,T_*).
\end{equation}
If $\xi$ is a unit exponential random variable, then
\[
 S\equaldis\sigma_*^2\xi,
\]
and $S$ is independent of $T_*$.  More explicitly, for $i\in I$ we have 
\begin{equation}
\label{eq:location-weights}
 \pk{T_*=i}=:\pi_i
 =
 \frac{\displaystyle
 \int_{\substack{\vk z_I>\vk0_I\\ z_i<z_j,\ j\in I\setminus\{i\}}}
 \psi_I(\vk z_I)e^{-\vk\theta^\top\vk z_I}\,d\vk z_I}
 {\displaystyle
 \int_{\vk z_I>\vk0_I}
 \psi_I(\vk z_I)e^{-\vk\theta^\top\vk z_I}\,d\vk z_I}.
\end{equation}
If $L_I^0=\emptyset$, the coordinates $(Z_I)_i$, $i\in I$, are
independent exponentials with rates $\theta_i$, and
\[
 \pi_i=\frac{\theta_i}{\sum_{j\in I}\theta_j}
      =\frac{\theta_i}{\mu^2}=p_i.
\]
\end{korr}

\begin{remark}\label{rem:dirichlet-location-weights}
The exponential overshoot in \nekorr{korr:minimum-location} is
universal within the present Dirichlet model, but the location weights
need not equal the quadratic-programming weights $p_i$.  When
$L_I^0\ne\emptyset$, the factor $\psi_I$ in
\eqref{eq:location-weights} may tilt the minimiser location.  Thus the
general location law depends on the Dirichlet parameters and on the
linear transformation $A$, not only on the shape matrix
$\Sigma=A^\top A$.  Accordingly, $\sigma_*^2$ is a shape-energy
quantity here; its interpretation as covariance energy is specific to
the Gaussian case.  When $L_I^0=\emptyset$, only the angular tilt of
the location law disappears, and $\pi_i=p_i$.
\end{remark}

\begin{remark}[Exact Gaussian recovery]\label{rem:gaussian-minimum}
Let $\X\sim N_k(\vk{0},\Sigma)$ be represented as in
Remark~\ref{rem:gaussian}, so that $\alpha_i=1/2$, $1\le i\le k$,
and $R^2\sim\chi_k^2$.  Put
\[
 \Gamma_P:=\Sigma_{PP}
 -\Sigma_{PI}\SIIIM\Sigma_{IP},\qquad
 c_P:=\pk{\vk G_P>\vk0_P},
\]
where $\vk G_P$ is centred Gaussian with covariance $\Gamma_P$;
set $c_P=1$ when $P=\emptyset$.  Then $L=\emptyset$,
$w(r)=r$, $c_n=u_n$, $d_n=1$, and
\[
 \bs_J=\Sigma_{JI}\SIIIM\vk{1}_I.
\]
When $J\ne\emptyset$, block inversion gives
\[
 \bigl((\Sigma^{-1})_{JJ}\bigr)^{-1}
 =\Sigma_{JJ}-\Sigma_{JI}\SIIIM\Sigma_{IJ}.
\]
Hence the $P$-subvector of the residual limit in
\eqref{eq:residual-limit-density} has the law of $\vk G_P$ conditioned
to be positive.
Consequently, \eqref{eq:conditional-active-passive} becomes
\[
 \mathcal L\left(
 \left(
 (u_n(X_i-u_n))_{i\in I},
 (X_j-u_n)_{j\in P},
 \frac{\X_N}{u_n}
 \right)
 \,\middle|\,\underline X>u_n\right)
 \Longrightarrow
 \mathcal L\left(
 (\xi_i/\theta_i)_{i\in I},
 \vk G_P^+,
 \bs_N
 \right), \qquad n\to\infty,
\]
where the $\xi_i$ are independent unit exponentials,
$\vk G_P^+\equaldis\vk G_P\mid\vk G_P>\vk0_P$, and
$\vk G_P^+$ is independent of the exponentials.  Furthermore, as $n\to\infty$
\[
 \mathcal L\left(
 \left(u_n(\underline X-u_n),T\right)
 \,\middle|\,\underline X>u_n\right)
 \Longrightarrow
 \mathcal L(\sigma_*^2\xi,T_*),
\qquad
 \pk{T_*=i}=p_i,
\]
with independent limiting components.  Finally, we have as $n\to\infty$
\[
 \pk{\underline X>u_n}
 \sim
 \frac{c_P}
 {(2\pi)^{|I|/2}\det(\Sigma_{II})^{1/2}
  \prod_{i\in I}\theta_i}
 u_n^{-|I|}
 \exp\left\{-\frac{u_n^2}{2\sigma_*^2}\right\}.
\]
Thus \netheo{theo:main1} and
\nekorr{korr:minimum-location} reproduce the finite-dimensional
Gaussian tail asymptotic, conditional excess limit and
minimiser-location limit established in \cite{HN26}.
\end{remark}

\def\xomega{{\IF}}
\def\omegaA{{\IF}}
\section{Proofs}

\prooftheo{theo:main1} We give the details, including the upper bound.  Let $G$ be the distribution function supplied by \nelem{lem:auxG}, so that
\BQN\label{eq:def:as:fd}
\overline{G}(u) &=&(1+o(1))
uw(u)\overline{F}(u) \frac{\Gamma(\barAL)}{ 2\Gamma(\barAL+1)}, \quad u\to \infty
\EQN
and let $\vk{Y}^*\sim \kal{GSD}(k+1, \vk{\alpha}^*, G)$, where $\vk{\alpha}^*=(\alpha_1 \ldot \alpha_k,1)^\top$.  Put $K_0:=\njk$ and $\Y:=A^\top\vk{Y}^*_{K_0}$.  By \nelem{lem:dist}, $\Y$ has density
\BQNY
h(\x)&= & \prod_{i=1}^k  \abs{(C \x)_i}^{2 \alpha_i-1}
\frac{\Gamma(\barAL+1)}{\prod_{i=1}^k\Gamma(\alpha_i)  \abs{\SI}^{1/2}}
\int_{\normp{\x}}^\omegaA
 z^{-2 \barAL } \, d G(z),  \quad \x\in\R^k,
\EQNY
where $\normp{\x}:=(\x^\top \SIM \x)^{1/2}$.  The subvector representation in \nelem{lem:dist} gives
$\Y\equaldis A^\top(R_G\sqrt B)\U$, where $R_G\sim G$, $B\sim\kal{B}_{\barAL,1}$ and $\U\sim\LPSDA$ are independent.  By \nelem{lem:auxG},
$\pk{R_G\sqrt{B}>u}\sim\overline F(u)$.  Since $\thry$ is eventually componentwise positive and its active coordinates diverge, \nelem{lemA} yields
\BQN \label{eq:lastX}
\pk{\X > \thry}&=& (1+ o(1)) \pk{\Y > \thry}, \quad n\to \IF.
\EQN
It remains to evaluate the right-hand side.
Let $d_J:=\abs{J}$ and define $\vk v_n$ by
\BQNY
(\vk v_n)_I=\wtny\vk{1}_I,\qquad
(\vk v_n)_J=\left(\frac{\wtny}{\ntny}\right)^{1/2}\vk{1}_J
\EQNY
when $J$ is non-empty; if $J=\emptyset$ only the first relation is used.  With the change of variables
\[
\x=u_n\bs+\y/\vk v_n
\]
(componentwise division), the Jacobian is
\BQN\label{eq:jacobian}
d\x=\ntny^{d_J/2}\wtny^{-\abs{I}-d_J/2}\,d\y.
\EQN
Set
\BQN\label{eq:Sdef}
H_0(\y):=\vk{u}_I^\top\SIIIM\y_I+\frac12\y_J^\top(\SIM)_{JJ}\y_J,
\EQN
omitting the second term if $J=\emptyset$.  From \nelem{prop:pre:gaus},
$\bs_I=\b_I$, $\bs_J=\SI_{JI}\SIIIM\b_I$ when $J\ne\emptyset$, and
$\x^\top\SIM\bs=\x_I^\top\SIIIM\b_I$ for all $\x$.  Therefore, locally uniformly in $\y$,
\BQN\label{eq:localnorm}
\normp{u_n\bs+\y/\vk v_n}
&=&\ntny+\frac{H_0(\y)+o(1)}{\wtny},\quad n\to\IF.
\EQN
Write $\beta_i:=2\alpha_i-1$ and put
\[
 \vk q_n:=\vk v_n(\thry-u_n\bs),\qquad
 D_n:=\{\y:\y>\vk q_n\},\qquad
 D:=\{\y_I>\vk q_I,\ \y_J>\vk q_J\}.
\]
As usual, a component $(q_J)_j=-\infty$ imposes no restriction in the definition of $D$, and all $J$-terms are omitted when $J=\emptyset$.  Since $(\vk q_n)_I\to\vk q_I$, there is a constant $b_0<\infty$ such that, for all sufficiently large $n$,
\BQN\label{eq:active-lower-bound}
 y_i&\ge&-b_0,\qquad i\in I,\quad \y\in D_n.
\EQN

We now prove the uniform integrability needed below.  Write $\mu:=\normS{\b_I}$ and let
\[
 M^\circ:=\{i\in M:\beta_i\ne0\},
 \qquad \eta_i:=\frac{(C\bs)_i}{\mu},\quad i\in M^\circ.
\]
Notice that $\eta_i\ne0$ for every $i\in M^\circ$.  Define the normalised angular factor
\begin{align}
 \mathcal A_n(\y)
 &:={}
 \prod_{i\in M^\circ}
 \left|
 \frac{\eta_i+\lambda_n^{-1}C_{iI}\y_I+
                 \lambda_n^{-1/2}C_{iJ}\y_J}{\eta_i}
 \right|^{\beta_i}                                      \notag\\
 &\quad\times
 \prod_{i\in L_I^0}|C_{iI}\y_I|^{\beta_i}
 \prod_{i\in L_J^0}
 \left|C_{iJ}\y_J+\lambda_n^{-1/2}C_{iI}\y_I\right|^{\beta_i}.
 \label{eq:normalised-angular}
\end{align}
Then the following identity is exact:
\begin{align}
 \prod_{i=1}^k
 \left|\bigl(C(u_n\bs+\y/\vk v_n)\bigr)_i\right|^{\beta_i}
 &={}
 \tau_M^*\ntny^{\sum_{i\in M}\beta_i}
 \left(\frac{\ntny}{\wtny}\right)^{
             \frac12\sum_{i\in L_J^0}\beta_i}
 \wtny^{-\sum_{i\in L_I^0}\beta_i}
 \mathcal A_n(\y).
 \label{eq:angular-exact-factorisation}
\end{align}
Also put
\[
 \mathcal K_n(\y):=
 \ntny^{2\barAL}
 \int_{\normp{u_n\bs+\y/\vk v_n}}^\omegaA
 z^{-2\barAL}\frac{dG(z)}{\overline G(\ntny)} .
\]
We first record the weighted-tail equivalence
\begin{equation}\label{eq:weighted-G-tail}
 \int_r^\IF z^{-p_0}\,dG(z)
 \sim r^{-p_0}\overline G(r),\qquad r\to\IF,
 \quad p_0:=2\barAL .
\end{equation}
Indeed, for every $a>1$, continuity of the $G$ supplied by Lemma
\ref{lem:auxG} and monotonicity give
\[
 a^{-p_0}r^{-p_0}\{\overline G(r)-\overline G(ar)\}
 \le \int_r^\IF z^{-p_0}\,dG(z)
 \le r^{-p_0}\overline G(r).
\]
Since $G\in GMDA(w)$ has a rapidly varying tail, first letting
$r\to\infty$ and then $a\downarrow1$ proves
\eqref{eq:weighted-G-tail}.  Applying this equivalence with
$r=\normp{u_n\bs+\y/\vk v_n}$, and then using
\eqref{eq:localnorm}, $\normp{u_n\bs+\y/\vk v_n}/\ntny\to1$, and the
locally uniform Gumbel convergence for $G$, gives, locally uniformly in
$\y$,
\BQN\label{eq:innerlimit}
 \mathcal K_n(\y)&\longrightarrow&\exp(-H_0(\y)), \qquad n\to \IF.
\EQN

We first record a global bound for $\mathcal K_n$.  Set
\[
 \rho_n(\y):=\normp{u_n\bs+\y/\vk v_n},\qquad
 \Phi(\y):=1+\sum_{i\in I}y_i^+ +\norm{\y_J}^2.
\]
By \eqref{eq:new}, with $\ell:=\SIIIM\vk u_I>\vk0_I$, the exact quadratic expansion is
\begin{align}
 \frac{\wtny}{\ntny}\bigl(\rho_n(\y)^2-\ntny^2\bigr)
 &={}
 2\ell^\top\y_I+
 \frac{\wtny}{\ntny}
 \left(\frac{\y}{\vk v_n}\right)^\top
 \SIM
 \left(\frac{\y}{\vk v_n}\right).
 \label{eq:global-quadratic}
\end{align}
Positive definiteness of $\SIM$, \eqref{eq:active-lower-bound}, and $\lambda_n=\ntny\wtny\to\infty$ therefore imply, uniformly for $\y\in D_n$
\BQN\label{eq:global-coercivity}
 \frac{\wtny}{\ntny}\bigl(\rho_n(\y)^2-\ntny^2\bigr)
 &\ge& c_0\Phi(\y)-C_0
\EQN
for some constants $c_0,C_0>0$.  Indeed, the quadratic term in \eqref{eq:global-quadratic} is bounded below by
\[
 c\left(\frac{\norm{\y_I}^2}{\lambda_n}+\norm{\y_J}^2\right),
\]
whereas $\ell^\top\y_I$, on the set \eqref{eq:active-lower-bound}, controls $\sum_{i\in I}y_i^+$ up to an additive constant.

For $\Phi(\y)$ sufficiently large, \eqref{eq:global-coercivity} gives $\rho_n(\y)\ge\ntny$.  If
\[
 \Delta_n(\y):=\frac{\wtny}{\ntny}
 \bigl(\rho_n(\y)^2-\ntny^2\bigr),
\]
then
\[
 \wtny\bigl(\rho_n(\y)-\ntny\bigr)
 =\lambda_n\left(\sqrt{1+\frac{\Delta_n(\y)}{\lambda_n}}-1\right).
\]
Since $\lambda_n\to\infty$, \eqref{eq:global-coercivity} consequently implies
\BQN\label{eq:global-radial-shift}
 \wtny\bigl(\rho_n(\y)-\ntny\bigr)
 &\ge& c_1\sqrt{\Phi(\y)}
\EQN
whenever $\Phi(\y)$ is sufficiently large.  For bounded $\Phi(\y)$,
\eqref{eq:global-coercivity} instead gives
$\rho_n(\y)\ge\ntny-C/\wtny$.  Since $\ntny\wtny\to\infty$,
\[
 \mathcal K_n(\y)
 \le
 \left(\frac{\ntny}{\rho_n(\y)}\right)^{2\barAL}
 \frac{\overline G(\rho_n(\y))}{\overline G(\ntny)}
 \le C .
\]
Here the last bound follows by monotonicity when
$\rho_n(\y)\ge\ntny$, and by the locally uniform Gumbel convergence at
shifts in $[-C,0]$ otherwise.  When $\Phi(\y)$ is sufficiently large,
$\rho_n(\y)\ge\ntny$ and
$\wtny(\rho_n(\y)-\ntny)\ge c_1\sqrt{\Phi(\y)}$.
The one-sided bound \eqref{eq:potter-gmda}, applied to $G$ with a
sufficiently small Potter parameter, therefore shows that for every
$R>0$ there is a constant $K_R$ such that
\BQN\label{eq:radial-global-bound}
 \mathcal K_n(\y)&\le&K_R\Phi(\y)^{-R},
 \qquad \y\in D_n,
\EQN
uniformly for all sufficiently large $n$.

We next control all singular linear factors.  We use the following elementary fact.  If $r\ge1$, $B$ is an $r\times k$ matrix of full row rank, $\gamma_1\ldot\gamma_r\in(-1,0)$, and $Q$ is a unit cube in $\R^k$, then
\BQN\label{eq:affine-power-bound}
 \sup_{\vk d\in\R^r}
 \int_Q\prod_{j=1}^r|d_j+(B\y)_j|^{\gamma_j}\,d\y
 &\le&C\{\sigma_{\min}(B)\}^{-r},
\EQN
provided $\norm{B}$ is bounded above.  To see this, complete the rows of $B$ by an orthonormal basis of their orthogonal complement to obtain a nonsingular $k\times k$ matrix $\widetilde B$.  Then
\[
 |\det\widetilde B|=\{\det(BB^\top)\}^{1/2}
 \ge\{\sigma_{\min}(B)\}^{r}.
\]
After the change of variables $\vk z=\widetilde B\y$, the image of $Q$ has uniformly bounded diameter.  Moreover,
\[
 \sup_{a\in\R}\int_{a-L}^{a+L}|z|^\gamma\,dz<\infty,
 \qquad -1<\gamma<0.
\]
This proves \eqref{eq:affine-power-bound}.

Let
\[
 H:=\{i:\beta_i<0\},\qquad r:=|H|,\qquad
 d_+:=\sum_{i=1}^k\max(\beta_i,0).
\]
If $H=\emptyset$, the regular-cube estimate below applies to every cube and the irregular-cube argument is void.  Thus the latter argument is understood only when $r\ge1$.
Partition $\R^k$ into unit lattice cubes.  For a cube $Q$ meeting $D_n$, put
\[
 \Phi_Q:=\inf_{\y\in Q\cap D_n}\Phi(\y).
\]
The values of $\Phi$ on $Q\cap D_n$ are comparable, up to constants independent of $Q$ and $n$, with $\Phi_Q$.

Call $Q$ regular if, for every $i\in H\cap M^\circ$
\[
 \left|\eta_i+\lambda_n^{-1}C_{iI}\y_I+
                 \lambda_n^{-1/2}C_{iJ}\y_J\right|
 \ge\frac{|\eta_i|}{2},\qquad \y\in Q.
\]
On a regular cube all negative-power $M^\circ$-factors are uniformly bounded.  For fixed $\y_I$, \eqref{eq:affine-power-bound}, applied to $C_{L_J^-,J}$, controls the negative-power $L_J^0$-factors uniformly in the affine shift.  A second application to $C_{L_I^-,I}$ controls the negative-power $L_I^0$-factors.  The required full-row-rank properties are precisely the rank assumptions of the theorem.  All non-negative powers are bounded on $Q$ by $C\Phi_Q^{d_+}$.  Hence
\BQN\label{eq:regular-cube-bound}
 \int_{Q\cap D_n}\mathcal A_n(\y)\,d\y
 &\le&C\Phi_Q^{d_+}
\EQN
on every regular cube.

Suppose now that $Q$ is not regular.  Then, for some $i\in H\cap M^\circ$ and some $\y\in Q$,
\[
 \lambda_n^{-1}|C_{iI}\y_I|
 +\lambda_n^{-1/2}|C_{iJ}\y_J|
 \ge\frac{|\eta_i|}{2}.
\]
Since $Q$ meets $D_n$ and its diameter is bounded, \eqref{eq:active-lower-bound} implies
\BQN\label{eq:irregular-cube-far}
 \Phi_Q&\ge&c\lambda_n.
\EQN

It remains to control the simultaneous negative-power singularities on such a cube.  Up to fixed non-zero row multipliers, the coefficient matrix of all affine forms indexed by $H$ is
\[
 B_n=\mathcal R_n C_H
 \operatorname{diag}\left(
    \lambda_n^{-1}I_{|I|},\lambda_n^{-1/2}I_{|J|}
 \right),
\]
where the diagonal entry of $\mathcal R_n$ is $1$, $\lambda_n^{1/2}$, or $\lambda_n$ according as the corresponding row belongs to $M^\circ$, $L_J^0$, or $L_I^0$.  The rows of $C_H$ are linearly independent because $C$ is nonsingular.  Choose an $r\times r$ non-zero minor $C_{H,K}$.  The corresponding minor of $B_n$ satisfies
\[
 \bigl|\det\bigl((B_n)_{:,K}\bigr)\bigr|
 \ge c\lambda_n^{-r}.
\]
All entries of $B_n$ are uniformly bounded.  By the Cauchy--Binet
formula,
\[
 \det(B_nB_n^\top)
 =\sum_{\substack{K'\subset\{1,\ldots,k\}\\ \abs{K'}=r}}
 \det\bigl((B_n)_{:,K'}\bigr)^2
 \ge c\lambda_n^{-2r}.
\]
The remaining $r-1$ singular values are uniformly bounded, and hence
$\sigma_{\min}(B_n)\ge c\lambda_n^{-r}$.  Applying
\eqref{eq:affine-power-bound} to all negative-power factors, and then
using \eqref{eq:irregular-cube-far}, gives
\[
 \int_{Q\cap D_n}\mathcal A_n(\y)\,d\y
 \le C\lambda_n^{r^2}\Phi_Q^{d_+}
 \le C\Phi_Q^{r^2+d_+}.
\]
Together with \eqref{eq:regular-cube-bound}, this proves that, with $d_0:=k^2+d_+$,
\BQN\label{eq:all-cube-bound}
 \int_{Q\cap D_n}\mathcal A_n(\y)\,d\y
 &\le&C\Phi_Q^{d_0}
\EQN
for every unit cube $Q$ meeting $D_n$, uniformly in $n$.

Choose $R>d_0+k+2$ in \eqref{eq:radial-global-bound}. Combining \eqref{eq:radial-global-bound} and \eqref{eq:all-cube-bound} yields
\[
 \int_{Q\cap D_n}\mathcal A_n(\y)\mathcal K_n(\y)\,d\y
 \le C\Phi_Q^{d_0-R}.
\]
The series of the right-hand side over all unit cubes satisfying $y_i\ge-b_0-1$, $i\in I$, is finite.  Indeed, the number of such cubes with $\sum_{i\in I}y_i^++\norm{\y_J}^2\le s$ is $O(s^{|I|+|J|/2})$.  We have therefore proved
\BQN\label{eq:uniform-tail-negligible}
 \lim_{R_0\to\infty}\ \sup_{n\ge n_0}
 \int_{D_n\cap\{\norm{\y}>R_0\}}
 \mathcal A_n(\y)\mathcal K_n(\y)\,d\y&=&0.
\EQN

Finally, on every fixed ball $\{\norm{\y}\le R_0\}$,
$\boldsymbol 1_{D_n}\to\boldsymbol 1_D$ almost everywhere.  To make
the singular-factor step explicit, when $L_J^-\ne\emptyset$ choose a
right inverse $T_J$ of $C_{L_J^-,J}$.  The affine shift
$\lambda_n^{-1/2}C_{L_J^-,I}\y_I$ tends to zero uniformly for bounded
$\y_I$, and the change
$\y_J\mapsto\y_J+T_J\lambda_n^{-1/2}C_{L_J^-,I}\y_I$ reduces the
negative-power $J$-factors to ordinary translations of an
$L^1_{\mathrm{loc}}$ function.  Translation continuity, followed by
Fubini with the locally integrable $L_I^-$-factors, proves the required
joint convergence.  The negative-power $M^\circ$-factors converge
locally uniformly because the corresponding $\eta_i$ are non-zero.
Together with the local-uniform convergence in \eqref{eq:innerlimit}
and local boundedness of the non-negative-power factors, this gives
\[
 \boldsymbol 1_{D_n}\mathcal A_n\mathcal K_n
 \longrightarrow
 \boldsymbol 1_D
 \prod_{i\in L_I^0}|C_{iI}\y_I|^{\beta_i}
 \prod_{i\in L_J^0}|C_{iJ}\y_J|^{\beta_i}
 e^{-H_0(\y)}, \qquad n\to \IF
\]
in $L^1$ on that ball.  Letting first $n\to\infty$ and then $R_0\to\infty$ in view of \eqref{eq:uniform-tail-negligible} proves
\[
 \int_{D_n}\mathcal A_n(\y)\mathcal K_n(\y)\,d\y
 \longrightarrow\tau_{L_I^0,L_J^0}, \qquad n\to \IF.
\]
The same rank changes of variables give local integrability of every
negative-power factor.  At infinity,
$\SIIIM\vk{u}_I>\vk{0}_I$ and the positive definiteness of
$(\SIM)_{JJ}$ control the $I$- and $J$-directions, respectively,
including coordinates for which $(q_J)_j=-\infty$; the non-negative
powers grow only polynomially.  Thus
$\tau_{L_I^0,L_J^0}<\infty$.  Its domain has non-empty interior and
the relevant row forms vanish only on finitely many proper
hyperplanes, so $\tau_{L_I^0,L_J^0}>0$.
Returning the deterministic factors in \eqref{eq:angular-exact-factorisation} gives the desired convergence. Consequently,
\BQNY
\lefteqn{\int_{\y>\vk v_n(\thry-u_n\bs)}
\prod_{i=1}^k\abs{\bigl(C(u_n\bs+\y/\vk v_n)\bigr)_i}^{\beta_i}
\int_{\normp{u_n\bs+\y/\vk v_n}}^\omegaA z^{-2\barAL}\frac{dG(z)}{\overline G(\ntny)}\,d\y}\\
&=&(1+o(1))\tau_M^*\ntny^{-2\barAL+\sum_{i\in M}\beta_i+\frac12\sum_{i\in L_J^0}\beta_i}
\wtny^{-\frac12\sum_{i\in L_J^0}\beta_i-\sum_{i\in L_I^0}\beta_i}
\tau_{L_I^0,L_J^0}.
\EQNY
Combining this with \eqref{eq:jacobian} and the density constant gives
\BQNY
\pk{\Y>\thry}
&=&(1+o(1))\frac{\Gamma(\barAL+1)}{\prod_{i=1}^k\Gamma(\alpha_i)\abs{\SI}^{1/2}}
\overline G(\ntny)\tau_M^*\tau_{L_I^0,L_J^0}\\
&&\times
\ntny^{d_J/2-2\barAL+\sum_{i\in M}\beta_i+\frac12\sum_{i\in L_J^0}\beta_i}
\wtny^{-\abs{I}-d_J/2-\frac12\sum_{i\in L_J^0}\beta_i-\sum_{i\in L_I^0}\beta_i}.
\EQNY
Finally use \eqref{eq:def:as:fd}.  Since
$\sum_{i\in E}\beta_i=2\overline{\alpha}_E-\abs{E}$ for every
$E\subset\njk$, the exponent of $\ntny$ and the exponent of $\wtny$
after substituting \eqref{eq:def:as:fd} are both equal to $\kappa$ in
\eqref{eq:kappa} implying
\BQNY
\pk{\Y>\thry}
&=&(1+o(1))\tau_M^*\tau_{L_I^0,L_J^0}
\frac{\Gamma(\barAL)}{2\prod_{i=1}^k\Gamma(\alpha_i)\abs{\SI}^{1/2}}
(\ntny\wtny)^{\kappa}\overline F(\ntny),
\EQNY
and \eqref{eq:lastX} proves \eqref{eq:theo1:gen}.  The two displayed special cases follow by factorising the integral $\tau_{L_I^0,L_J^0}$ when $L_I^0=\emptyset$, and by putting $J=\emptyset$, respectively. \QED

\bigskip
\proofkorr{korr:ell-comparison}
Apply \netheo{theo:main1} to $\X$ and to $\X^*$ and take the ratio.
When $L=\emptyset$, their local integrals and powers of $\lambda_n$
coincide. \QED

\bigskip
\proofkorr{korr:excess}
For $\x\in[0,\infty)^k$ 
\[
 \pk{\wny\X[\thry]>\x}
 =\frac{\pk{\X>\thry+\x/\wny}}{\pk{\X>\thry}}.
\]
The threshold $\thry+\x/\wny$ satisfies the assumptions of
\netheo{theo:main1} with $\vk q$ replaced by $\vk q+\x$.
The factors $\tau_M^*$, the power of $\lambda_n$, and
$\overline F(\ntny)$ are the same in numerator and denominator, so
they cancel.  This proves \eqref{eq:excesslimit}.  If
$L_I^0=\emptyset$, the integral in the $I$-coordinates is elementary
and yields \eqref{eq:excessfactor}.  The same argument applies when
some components of $\vk q_J$ equal $-\infty$, with the convention
stated after the corollary. \QED

\bigskip
\proofkorr{korr:minimum-location}
The assertions concerning $\vk{p}$ follow either from the duality
between the two quadratic programmes or directly from
\[
 (\Sigma\vk{p})_I=\sigma_*^2\vk{1}_I,\qquad
 (\Sigma\vk{p})_J=\sigma_*^2\bs_J\ge\sigma_*^2\vk{1}_J.
\]
Indeed, these are the optimality conditions for the strictly convex
problem defining $\vk{p}$, and
$\vk{p}^\top\Sigma\vk{p}=\sigma_*^2$.

We next prove the conditional rescaled-vector limit.  Put
\[
 \vk{q}^0_I=\vk{0}_I,\qquad
 \vk{q}^0_P=\vk{0}_P,\qquad
 \vk{q}^0_N=-\vk{\infty}_N,
\]
where the last two subvectors are omitted when the corresponding sets
are empty.  For fixed $\vk{x}_I\ge\vk{0}_I$, $\vk{y}_P\ge\vk{0}_P$ and
$\vk{y}_N\in\R^{|N|}$, define a threshold $\vk{t}_n(\vk{x}_I,\vk{y}_J)$
by
\[
 \bigl(\vk{t}_n(\vk{x}_I,\vk{y}_J)\bigr)_I
 =u_n\vk{1}_I+\frac{\vk{x}_I}{c_n},\qquad
 \bigl(\vk{t}_n(\vk{x}_I,\vk{y}_J)\bigr)_J
 =u_n\bs_J+\frac{\vk{y}_J}{d_n}.
\]
Since $r_nw(r_n)\to\infty$ as $n\to\infty$, the normalising constants   satisfy
\[
 u_nd_n=\frac{\sqrt{r_nw(r_n)}}{\mu}\longrightarrow\infty, \qquad n\to \IF.
\]
Consequently, for all sufficiently large $n$ the event
$\{\X>\vk{t}_n(\vk{x}_I,\vk{y}_J)\}$ is contained in
$\{\underline X>u_n\}$.  The active and residual limits of this
threshold in \netheo{theo:main1} are
\[
 w(r_n)\bigl(\vk{t}_n-u_n\bs\bigr)_I=\mu\vk{x}_I,
 \qquad
 d_n\bigl(\vk{t}_n-u_n\bs\bigr)_J=\vk{y}_J.
\]
For the denominator threshold $u_n\vk{1}$, they are
$\vk{q}^0_I,\vk{q}^0_P,\vk{q}^0_N$.  Taking the ratio of the two
expansions in \eqref{eq:theo1:gen}, and then making the change of
variables $\vk{y}_I=\mu\vk{z}_I$ in the active integral, gives
\begin{align*}
 &\pk{\,
 c_n(\X_I-u_n\vk{1}_I)>\vk{x}_I,\ 
 d_n(\X_J-u_n\bs_J)>\vk{y}_J
 \mid \underline X>u_n\,}\\
 &\hspace*{1.2cm}\longrightarrow
 \pk{\vk{Z}_I>\vk{x}_I,\ \vk{Z}_J^{(P)}>\vk{y}_J}, \qquad n\to \IF.
\end{align*}
Here the inequalities in the $N$-coordinates are unrestricted in the
denominator.  The limiting law has the product density
\eqref{eq:active-limit-density}--\eqref{eq:residual-limit-density}.
In particular, its two blocks are independent, and its normalising
constants are finite and positive by \netheo{theo:main1}.  Since this
law is absolutely continuous, convergence on these upper orthants
implies \eqref{eq:conditional-local-min}.

The $N$-part of \eqref{eq:conditional-local-min} is tight and
$u_nd_n\to\infty$, hence
\[
 \frac{\X_N}{u_n}
 =\bs_N+\frac{d_n(\X_N-u_n\bs_N)}{u_nd_n}
 \longrightarrow\bs_N
\]
in conditional probability, which  proves
\eqref{eq:conditional-active-passive}.  Moreover, we have as $n\to\infty$
\[
 \frac{c_n}{d_n}=\frac{\sqrt{r_nw(r_n)}}{\mu}\longrightarrow\infty.
\]
Since $(\vk{Z}_J^{(P)})_P>\vk{0}_P$ almost surely 
\[
 c_n(\X_P-u_n\vk{1}_P)
 =\frac{c_n}{d_n}\,d_n(\X_P-u_n\vk{1}_P)
 \longrightarrow+\vk{\infty}_P
\]
in conditional probability and further  
\[
 c_nu_n=\frac{r_nw(r_n)}{\mu^2}\longrightarrow\infty,\qquad
 \frac{\X_N}{u_n}-\vk{1}_N
 \longrightarrow\bs_N-\vk{1}_N>\vk{0}_N, \qquad n\to\infty
\]
and thus  $c_n(\X_N-u_n\vk{1}_N)\to+\vk{\infty}_N$ in conditional
probability.  Thus, with conditional probability
tending to one, both $\underline X$ and $T$ are determined by the
active coordinates.  The density \eqref{eq:active-limit-density} is
absolutely continuous, so its minimum is unique almost surely and the
continuous-mapping theorem gives \eqref{eq:minimum-location}.  The
original minimum is also almost surely unique: for $i\ne j$,
$X_i-X_j$ is a non-zero linear functional of $R\U$, whose zero set has
zero symmetrised Dirichlet angular measure.

It remains to identify the law of $S$.  For every $\ell\in L_I^0$,
$C_{\ell J}=\vk{0}_J$ by definition, and hence
\[
 C_{\ell I}\vk{1}_I
 =C_{\ell I}\bs_I+C_{\ell J}\bs_J
 =(C\bs)_\ell=0.
\]
Therefore
\[
 \psi_I(\vk{z}_I+x\vk{1}_I)=\psi_I(\vk{z}_I),
 \qquad x\in\R,
\]
whenever both sides are defined, while
$\vk{1}_I^\top\vk{\theta}=\mu^2$.  Let
\[
 D_i:=\{\vk{z}_I>\vk{0}_I:z_i<z_j,\ 
                  j\in I\setminus\{i\}\}.
\]
For $x\ge0$, translation by $x\vk{1}_I$ in the integral over $D_i$
gives
\[
 \pk{S>x,T_*=i}
 =e^{-\mu^2x}\pk{T_*=i}
 =e^{-x/\sigma_*^2}\pi_i.
\]
This proves simultaneously that
$S\equaldis\sigma_*^2\xi$, that $S$ and $T_*$ are independent, and
that the weights are given by \eqref{eq:location-weights}.  Finally,
if $L_I^0=\emptyset$, \eqref{eq:active-limit-density} is the product
of exponential densities with rates $\theta_i$, so
$\pi_i=\theta_i/\mu^2=p_i$. \QED

\appendix
\section{Auxiliary results}\label{sec:auxiliary}  
Proof of \eqref{eq:resn}. A direct proof is given in Lemma 4.2 of \cite{EH2009a}. Another proof follows from the arguments of Lemma 2.2(b) in \cite{resnick2}. We record the short argument, since the same bound is used in the proof of the theorem.  For every sufficiently small $\ve\in(0,1)$ there are constants $c>0$ and $x_0$ such that
\BQN\label{eq:potter-gmda}
 \frac{\overline{F}(x + s/w(x))}{\overline{F}(x)} &\le & \frac{c}{(1 + \ve s)^{ \ve ^{-1}}},\qquad s\ge0,\ x\ge x_0 .
 \EQN
Let $r>1$ and put $s_x=(r-1)xw(x)$.  By \eqref{eq:uv}, $s_x\to\infty$.  Hence, for $x$ large,
\[
\frac{(xw(x))^\eta\overline F(rx)}{\overline F(x)}
\le c(xw(x))^\eta\{1+\ve(r-1)xw(x)\}^{-1/\ve}.
\]
If $\eta\le0$ this tends to zero immediately; if $\eta>0$, choose $\ve$ so small that $1/\ve>\eta$.  This proves \eqref{eq:resn}.  The same bound is valid for any distribution function in $GMDA(w)$, in particular for the auxiliary distribution function $G$ introduced below.
 
The following elementary construction supplies the radial distribution used in the proof.
\BL\label{lem:auxG} Suppose that Assumption A1 holds.  There exists a
continuous distribution function $G\in GMDA(w)$ such that
\eqref{eq:def:as:fd} holds.  If $R_G\sim G$ is independent of
$B\sim\kal{B}_{\barAL,1}$, then
\BQN\label{eq:beta_tail_equiv}
\pk{R_G\sqrt B>u}&\sim&\overline F(u),\quad u\to\IF .
\EQN
\EL

\prooflem{lem:auxG} Put
\[
q_u:=uw(u),\qquad
H(u):=q_u\overline F(u)
\frac{\Gamma(\barAL)}{2\Gamma(\barAL+1)},\quad u>0.
\]
The locally uniform convergence in \eqref{eq:rdfd} and the
self-neglecting property of $1/w$ imply
\[
 \frac{H(u+x/w(u))}{H(u)}\longrightarrow e^{-x}, \qquad u\to\infty      
\]
locally uniformly for $x\in\R$, and $H(u)\to0$.  Thus $H$ is
Gamma-varying.  By the standard monotone-equivalence theorem for
Gamma-varying functions, see \cite{dehaan,resnick2}, it has an
equivalent continuous non-increasing version on a sufficiently large
half-line.  Use this version as $\overline G$ there and extend it to a
continuous survival function on $[0,\infty)$ with
$\overline G(0)=1$.  Then $G(0)=0$ and as $u\to\infty$
\[
 \overline G(u)\sim H(u),
\]
and equivalence at both $u$ and $u+x/w(u)$ shows that
$G\in GMDA(w)$.  This proves \eqref{eq:def:as:fd}.

It remains to verify \eqref{eq:beta_tail_equiv}.  Since $B$ has density
$\barAL b^{\barAL-1}$ on $(0,1)$,
\[
 \pk{R_G\sqrt B>u}
 =\barAL\int_0^1\overline G(u/\sqrt b)b^{\barAL-1}\,db.
\]
Fix $\epsilon\in(0,1)$ and put
$r_\epsilon:=(1-\epsilon)^{-1/2}>1$.  The contribution of
$b\le1-\epsilon$ is at most
\[
 (1-\epsilon)^{\barAL}\overline G(r_\epsilon u)
 =o\left(\frac{\overline G(u)}{q_u}\right)
 =o(\overline F(u)),
\]
by \eqref{eq:resn} applied to $G$.  On $b\in(1-\epsilon,1)$ put
$s=q_u(b^{-1/2}-1)$ and
$S_u:=q_u\{(1-\epsilon)^{-1/2}-1\}$.  The change of variables gives
\[
 \barAL\int_{1-\epsilon}^1
 \overline G(u/\sqrt b)b^{\barAL-1}\,db
 =
 \frac{2\barAL}{q_u}\overline G(u)
 \int_0^{S_u}
 \frac{\overline G(u+s/w(u))}{\overline G(u)}
 \left(1+\frac{s}{q_u}\right)^{-(2\barAL+1)}\,ds .
\]
Here $S_u\to\infty$.  Extend the integrand by zero beyond $S_u$.
The local-uniform Gumbel convergence gives pointwise convergence to
$e^{-s}$, while \eqref{eq:potter-gmda}, with a Potter parameter
$\delta\in(0,1)$, gives the integrable bound
$C(1+\delta s)^{-1/\delta}$.  Dominated convergence therefore yields
\[
 \pk{R_G\sqrt B>u}
 \sim\frac{2\barAL}{q_u}\overline G(u)
 \sim\overline F(u), \qquad u\to\infty
\]
because $\Gamma(\barAL+1)=\barAL\Gamma(\barAL)$. \QED
 
The next lemma transfers equivalent radial tails through positive orthant events.
\BL\label{lemA} Let $R,R^*$ be positive random variables with infinite
upper endpoints such that
\BQN \limit{u} \frac{\pk{R>u}}{\pk{R^*>u}}&=& c \in (0,\IF).
\EQN
Let $\vk Z$ be a bounded random vector, independent of both radial variables.  If $\vk t_n\in(0,\infty)^k$ and $\max_i t_{n,i}\to\infty$, then, whenever the denominator is positive,
\BQN
 \frac{\pk{R\vk Z>\vk t_n}}{\pk{R^*\vk Z>\vk t_n}}&\longrightarrow& c.
\EQN
In particular, this applies to $\vk Z=A^\top\U$ with $A$ non-singular and $\U\sim\LPSDA$.
\EL

\prooflem{lemA} Choose $K\in(0,\infty)$ such that
$\max_i|Z_i|\le K$ almost surely.  If a component of $\vk Z$ is
non-positive, the event $R\vk Z>\vk t_n$ is empty.  On
$\{\vk Z>\vk0\}$ put
\[
 r_n(\vk Z):=\max_{1\le i\le k}\frac{t_{n,i}}{Z_i}.
\]
Then $r_n(\vk Z)\ge K^{-1}\max_i t_{n,i}\to\infty$ uniformly in $\vk Z$, and
\[
 \pk{R\vk Z>\vk t_n}
 =\E{\pk{R>r_n(\vk Z)}\,\boldsymbol 1_{\{\vk Z>\vk0\}}},
\]
with the analogous identity for $R^*$.  Since $a_n:=K^{-1}\max_i t_{n,i}\to\infty$, the assumed tail equivalence implies
\[
 \sup_{s\ge a_n}\left|\frac{\pk{R>s}}{\pk{R^*>s}}-c\right|\longrightarrow0.
\]
Hence the first integrand equals $(c+o(1))$ times the second, uniformly in $\vk Z$, and integration proves the claim. \QED
 
The next lemma can be found in Hashorva (2005).
 \BL\label{prop:pre:gaus}
Let $\Sigma \inr^{k \times k}$, $k\ge 2$, be positive definite and let $\bb\in \R^k \setminus (-\IF, 0]^k$ be a given vector.
Then the quadratic programming problem
$$ \qp{\SIM}{\bb}: \text{minimise $ \normS{\x}^2$ under the linear constraint } \x \ge \bb $$
has a unique solution $\bs$ defined by a   unique non-empty
index set $\IB\subset \njk$ such that
\BQN \label{eq:IJ}
\bs_{\IB}= \bb_{\IB},    \quad \SIIIM \bb_{\IB}&>&\vk{0}_I, \\
\label{eq:alfa} \min_{\x \ge  \bb} \normS{\x}^2=\min_{\x \ge  \bb}
\x^\top \SIM\x= \normS{\bs_I}^2
= \normS{{\bb_{\IB}}}^2 & = & \bb_{\IB}^\top \SIIIM\bb_{\IB}>0
\EQN
and if $\abs{I}< k$, then with $\JB :=\njk\setminus I$
\BQN
\label{eq:prop:A1} \bs_{\JB}&=&   - ((\SI^{-1})_{JJ})^{-1}(\SIM)_{JI}\bb_I=
\SIJI \SIIIM \bb_{\IB}\ge \bb_{\JB}.
\EQN
 Furthermore, for any $\x\inr^k$ we have
\BQN \label{eq:new}
\x^\top \SIM \bs= \x_I^\top
\SIIIM \bb_I= \x_I^\top \SIIIM\bs_I.
\EQN
\EL

\bigskip
 
The next result follows directly from the Dirichlet representation; see \cite{FangFang1990}.

\BL\label{lem:dist} Let $\X\sim \GSDAKF$ be a random vector in $\R^k,k\ge 2$ where
$F$ is a distribution function satisfying $F(0)=0$. For any non-empty index
$I\subset \njk$ with $m$ elements we have
\BQN \X_I \equaldis R_I \widetilde{\U}_I,
\EQN
with $R_I$ being independent of $\widetilde{\U}_I\sim
\kal{SD}(m, \vkal_I)$. If $m<k$, one may take
\[
 R_I=R\sqrt{B_I},\qquad
 B_I\sim\kal{B}_{\barALI,\barAL-\barALI},
\]
where $R$, $B_I$, and $\widetilde{\U}_I$ are mutually independent.  In particular,
$R_I$ possesses the density function $f_I$
defined by
\BQN \label{eq:density:ehen}
 f_I(z)&=& 2 z^{2\barALI-1} \frac{\Gamma(\barAL)}{\Gamma(\barALI)\Gamma(\barAL- \barALI)}
 \int_{z}^\IF (r^2- z^2 )^{\barAL- \barALI-1}r^{-2(\barAL- 1)} \, d F(r), \quad \forall z\in (0,\xomega).
\EQN
Furthermore, if $m<k$, then for any non-singular matrix $A\in
\R^{m \times m } $ the random vector $A^\top \X_I$ possesses the density function
$h_{A,I}$ given  by (set $C:=(A^\top)^{-1}, \SI:= A^\top A$, $\normp{\x_I}:=(\x_I^\top\SI^{-1}\x_I)^{1/2}$, and write $\abs{\SI}$ for the determinant of $\SI$)
\BQN\label{eq:hI}
h_{A,I}(\x_I)&= & \frac{\Gamma(\barAL)\prod_{i\in
I} \abs{(C \x_I)_i}^{2 \alpha_i-1}}{\prod_{i\in I}
\Gamma(\alpha_i)   \Gamma(\barAL- \barALI) \abs{\SI}^{1/2}}
 \int_{\normp{\x_I}}^\omegaA (r^2- \normp{\x_I}^2 )^{\barAL- \barALI-1}r^{-2(\barAL- 1)} \, d F(r),\notag \qquad\forall \x_I\inr^m.
\EQN
If $m=k$, then one may take $R_I=R$ and
$\widetilde{\U}_I=\U$, and the preceding density formulae involving
$\Gamma(\barAL-\barALI)$ are replaced by the original density representation of $\X$.
\EL

\end{document}